\documentclass[12pt]{amsart}

\def\@typesizes{%
       \or{5}{6.5}\or{6}{7.5}\or{7}{8.5}\or{8}{11}\or{9}{12}%
       \or{10}{13}
       \or{\@xipt}{14}\or{\@xiipt}{15}\or{\@xivpt}{18}%
       \or{\@xviipt}{20}\or{\@xxpt}{24}}

\usepackage{amsthm}
\usepackage{amsmath}
\usepackage{amssymb}
\usepackage{graphicx}
\usepackage{enumerate}
\usepackage{color}

\allowdisplaybreaks
\numberwithin{equation}{section}
\numberwithin{figure}{section}

\theoremstyle{plain}
\newtheorem{theorem}{ Theorem}[section]
\newtheorem{proposition}[theorem]{ Proposition}
\newtheorem{lemma}[theorem]{ Lemma}
\newtheorem{corollary}[theorem]{ Corollary}
\newtheorem{example}[theorem]{ Example}
\newtheorem{remark}[theorem]{ Remark}
\newtheorem{definition}[theorem]{ Definition}
\newtheorem{conjecture}{ Conjecture}

\newtheorem{lettertheorem}{ Theorem}

\numberwithin{equation}{section}

\def\BET{\begin{theorem}}
\def\ENT{\end{theorem}}
\def\BEP{\begin{proposition}}
\def\ENP{\end{proposition}}
\def\BEL{\begin{lemma}}
\def\ENL{\end{lemma}}
\def\BEC{\begin{corollary}}
\def\ENC{\end{corollary}}
\def\BEE{\begin{example} \rm}
\def\ENE{\end{example}}
\def\BER{\begin{remark} \rm}
\def\ENR{\end{remark}}
\def\BED{\begin{definition} \rm}
\def\END{\end{definition}}
\def\BECJ{\begin{conjecture}}
\def\ENCJ{\end{conjecture}}

\def\bea{\begin{eqnarray}}
\def\eea{\end{eqnarray}}

\def\beas{\begin{eqnarray*}}
\def\eeas{\end{eqnarray*}}

\def\beq{\begin{equation}}
\def\eeq{\end{equation}}

\def\beal{\begin{align*}}

\def\eeal{ \end{align*} }

\def\row{ \nonumber \\ & & }
\def\roweq{\nonumber \\ &=& }
\def\rowleq{\nonumber \\  & \leq & }
\def\rowgeq{\nonumber \\ & \geq & }

\def\bbC{{\mathbb C}}
\def\bbD{{\mathbb D}}

\def\bbN{{\mathbb N}}

\def\bbR{{\mathbb R}}

\def\bbZ{{\mathbb Z}}

\def\cP{{\mathcal P}}

\def\ef{\eqref}

\begin{document}

\title[Toeplitz operators]{A geometric condition for the invertibility of Toeplitz operators
on the Bergman space}

\author{$\check{\mbox{Z}}$eljko $\check{\mbox{C}}$u$\check{\mbox{C}}$kovi\'c}
\address{Department of Mathematics, University of Toledo, Toledo, OH 43606, USA}
\email{Zeljko.Cuckovic@utoledo.edu}

\author{Jari Taskinen}
\address{University of Helsinki, Department of Mathematics and Statistics,
P.O.Box 68, 00014 Helsinki, Finland}
\email{Jari.Taskinen@helsinki.fi}



\begin{abstract}
Invertibility of Toeplitz operators on the Bergman space and the related Douglas problem are long standing open problems.  In this paper we study the invertibility problem under the novel geometric condition on the image of the symbols, which relaxes the standard positivity condition.  We show that under our geometric assumption, the Toeplitz operator $T_\varphi$ is invertible if and only if the Berezin transform of $|\varphi|$ is invertible in $L^{\infty}$. It is well known that the Douglas problem is still open for harmonic functions.  We study a class of rather general harmonic polynomials and characterize the invertibility of the corresponding Toeplitz operators. We also give a number of related results and examples.
\end{abstract}

\maketitle

\section{Introduction.}
\label{sec1}

In the study of  Toeplitz operators $T_\varphi$ in the Bergman space, the characterization of the invertibility 
of $T_\varphi$ in terms of the properties of the symbol $\varphi$ and the related Douglas 
problem are important long standing open problems. Partial answers can be given by assuming the positivity of the 
symbol, which  is anyway a standard assumption in the literature used for a large number of results on the 
boundedness and other properties of Toeplitz operators in many types of analytic function spaces. 
In this paper we replace the positivity assumption by a novel, much weaker geometric condition 
and use it to study the invertibility of Toeplitz operators $T_\varphi$ in the Bergman space 
$A^2 $ on the open unit disc $\bbD$ of the complex plane. The condition concerns bounded symbols $\varphi$ 
on $\bbD$ and requires that Re\,$\varphi(z)$ is bounded from below by the 
non-negative  quantity $|{\rm Im} \varphi(z)|^2$. For such symbols,  we characterize
the invertibility of $T_\varphi$ in terms of the invertibility of the Berezin transform $B(|\varphi|)$ 
(Theorem \ref{th2.2}). As an application of the result, we consider a large class of harmonic polynomials
$P$ and characterize in Theorem \ref{th3.3} the invertibility of the operators $T_P$, thus  contributing to 
the Douglas problem, which is still open in the most important case of harmonic symbols. 
We will also give a number of other examples and variations of the result. 

In the sequel, we denote by $L^2 = L^2(\bbD)$ the Hilbert-space
of square integrable functions with respect to the normalized area measure $dA = 
\frac1\pi r dr d\theta$ on $\bbD$.  
We denote the inner produce of $A^2$ by $\langle f  ,g \rangle = \int_\bbD f \overline g dA$ 
and the norm by $\Vert f \Vert_2 = \sqrt{\langle f  , f \rangle}$,
where $f,g \in L^2$. The Bergman  space $A^2= A^2(\bbD)$ is the closed subspace of $L^2$ consisting of analytic functions. The orthogonal projection $P$ from
$L^2$ onto $A^2$ is called the Bergman projection, and it can be written 
as the integral operator
\beas
P f(z) = \int\limits_\bbD  K(z,w) f(w) dA(w), \ \ \ \mbox{where} \  K(z,w) = \frac{1}{(1- z \bar w)^2},
\ z \in \bbD, 
\eeas
is the Bergman kernel. The normalized kernel 
\beas
k_z(w) = \frac{1 - |z|^2}{(1- w \bar z)^2}
\eeas
has the property $\Vert k_z \Vert_2 = 1$ for all $z \in \bbD$. 

Let us denote by $L^\infty = L^\infty(\bbD)$ the Banach space of (essentially) bounded functions,
endowed with the sup-norm $\Vert \cdot \Vert_\infty$. We will consider Toeplitz operators $T_\varphi$
with  symbols $\varphi \in L^\infty$. The operator $T_\varphi$ is defined as the product $PM_\varphi$,
where $M_\varphi$ is the pointwise multiplier $M_\varphi: f \mapsto \varphi f$ so that we can write
\beas
T f(z) = \int\limits_\bbD \frac{\varphi(w) f(w)}{(1- z \bar w)^2} dA(w), \ \ \ f \in A^2, \ z \in \bbD. 
\eeas
Let $T$ be a bounded linear operator $A^2 \to A^2$ and let $g \in L^\infty$. The Berezin transform 
of $T$, as well as that of $g$, are complex-valued functions on
$\bbD \ni z$ defined by 
\bea
B(T)(z) = \langle T k_z , k_z \rangle \ \ \ \mbox{and} \ \ \
B(g)(z)= 
\int\limits_\bbD |k_z (w) |^2 g(w)  dA(w)    ,  \label{1.10}
\eea
respectively. It is known that $ |B(T)(z)| \leq \Vert T \Vert$ for all $z \in \bbD$, where $\Vert T\Vert$ denotes
the operator norm of $T$; see  Proposition 6.1. in \cite{Zh}.

R.~Douglas \cite{D} asked in the case of Hardy space $H^2$, if the invertibility of a Toeplitz operator 
with symbol belonging to $L^\infty (\partial \bbD)$ follows from  
the invertibility of the harmonic extension of the symbol from $\partial \bbD$ to 
$\bbD$. In the Bergman space setting the question thus reads, if the invertibility of $T_\varphi$ is 
implied by the invertibility of the Berezin transform of the symbol. We emphasize that the 
question is still open in the most important case that $\varphi$ is harmonic. However,  G.~McDonald and C.~Sundberg 
\cite{McS} showed that for a real harmonic symbol, the spectrum coincides with the interval 
$[\inf \varphi, \sup \varphi]$. Since $B( \varphi)$ and $\varphi$ coincide in this case, the 
result of \cite{McS} implies a positive answer to Douglas' question for real harmonic symbols. Moreover,
in the case of analytic symbols, the Toeplitz operator reduces to the 
pointwise multiplier, hence, the invertibility is equivalent with the existence of a bounded inverse $\varphi(z)^{-1}$ of the symbol   on the disc (see the argument in \cite{Hal}, p. 197).

The following sufficient and necessary conditions for the invertibility of $T_\varphi$ 
was proved by  D.~Luecking in \cite{Lue}. 
Here, we denote by $|A|$ the area of a set $A \subset \bbC$ and by $D(w,\varepsilon)$  the pseudohyperbolic disc with center $w$ and radius $\varepsilon$ (see \cite{Zh}, Section 4.2). 

\begin{lettertheorem}  \label{lem1.1}
For a non-negative $\varphi \in L^\infty$, the following conditions are equivalent.

\smallskip

\noindent $(i)$ The Toeplitz operator $T_\varphi$ is invertible on $A^2$.

\smallskip

\noindent $(ii)$ There exists a constant $\eta > 0$ such that 
\bea
\int\limits_\bbD |\varphi f|^2 dA \geq \eta \int\limits_\bbD |f|^2 dA \ \ \ \ \forall\, f \in A^2. 
\eea

\smallskip

\noindent $(iii)$  There exist $r,\delta > 0$ and $\varepsilon \in (0,1)$ such that
\bea
| G \cap D(w, \varepsilon) | > \delta |D(w, \varepsilon)| \ \ \ \ \forall\, w \in \bbD,
\eea
where $G = \{ z \in \bbD \, : \, \varphi(z) > r \}$.
\end{lettertheorem}

The formulation of Theorem \ref{lem1.1} actually coincides with Lemma 2.4. of \cite{ZZ16}. In \cite{ZZ16}, Corollary 3.5,  X.~Zhao and D.~Zheng showed that the answer to Douglas'  question is in general 
negative by constructing a counterexample with a symbol $\varphi$, which is a second order real polynomial 
of $|z|$. Also, the following result is Theorem 3.2 of \cite{ZZ16}.  

\begin{lettertheorem}  \label{thZZ3.2}
Let $\varphi \in L^\infty$ be non-negative. The operator $T_\varphi$ is invertible, if and only if the Berezin 
transform $B(\varphi)$ is invertible in $L^\infty$. 
\end{lettertheorem}

In Theorem 4.2 of \cite{ZZ16} it was proved that if $\varphi$ is a bounded harmonic function such that 
for some $\delta > 27/\sqrt{730} \approx 0.999315$ there holds
$ 
 \delta \Vert \varphi \Vert_\infty \leq |\varphi(z)| \leq  \Vert \varphi \Vert_\infty ,
$ 
for all $z \in \bbD$, then $T_\varphi$ is invertible. In \cite{ZZ19}, the authors improved this result: the 
constant $27/\sqrt{730}$ was replaced by $2 \sqrt{2} / 3$, which is smaller. 
The invertibility problem  for operators with harmonic symbols has also been 
considered in \cite{AC}, \cite{GZ}, \cite{GZZ} and \cite{Yo}. Moreover, 
in \cite{CucVas} the first named author and A.~Vasaturo characterized the invertibility of Toeplitz operators
$T_\mu$ where the symbol is a positive measure and showed that the invertibility of $T_\mu$ 
is equivalent to $\mu$ being a Reverse Carleson measure. This led to a necessary and sufficient 
condition for the 
invertibility of Toeplitz operators whose symbols are averaging functions of these Carleson measures.
The Douglas question for operators with measure symbols have recently also been considered in 
\cite{ZZZ}.

Here, our Theorem \ref{th2.2} 
generalizes Theorem \ref{thZZ3.2} by replacing the strong positivity assumption on the symbol by the
weaker geometric assumption; in addition, one needs to consider the Berezin transform $B(|\varphi|)$ instead
of $B(\varphi)$.  This result leads to 
a criterion for the invertibility of operators, the symbols of which belong to a rather 
general class of harmonic polynomials, see Theorem \ref{th3.3}.

In Section \ref{sec7} we present variations of the geometric condition, which are sufficient for the 
invertibility of $T_\varphi$. 
Section \ref{sec6} contains a reformulation of Theorem \ref{th2.2} for the case of harmonic
symbols and some more examples and remarks concering the Berezin transform. 
Also, in Theorem \ref{hX2} we prove that regular enough symbols, the  geometric condition implies that
the invertibility and Fredholmness properties of the operator$T_\varphi$ coincide.

\section{Preliminary results}  \label{sec2}

In this section we present the preliminary considerations needed for the proofs of the main results
in Sections \ref{sec3}--\ref{sec6}.

\BEL
\label{lem2.1}
Let $\varphi : \bbD \to \bbC$ be a bounded, measurable function such that $\Vert \varphi \Vert_\infty 
\leq 1$ and 
\bea
{\rm Re} \, \varphi(z) \geq |\varphi(z)|^2
\label{2.5d}
\eea
for all $z \in \bbD$. Moreover, assume that there exist 
$S> 0$, $M \geq 1$ such that the set  $G= \{  z \, : \, |\varphi (z)| > S \}$
has the property
\bea
\int\limits_\bbD |f|^2 dA < M \int\limits_G |f|^2 dA \ \ \
\forall \, f \in A^2.  \label{2.4}
\eea
Then, the Toeplitz-operator $T_\varphi$ is invertible on $A^2$, and there holds
the bound $\Vert T_\varphi^{-1} \Vert \leq MS^{-2}$. 
\ENL

Proof. Using \ef{2.5d} two times implies 
\bea
& & | 1 - \varphi(z)|^2 = 1- 2  {\rm Re}\, \varphi(z) + \varphi(z)^2
\leq  1- {\rm Re}\, \varphi(z) \leq 1- |\varphi(z)|^2 . \label{2.6a}
\eea
Following  the idea in the  proof of Corollary 3 of \cite{Lue}, the inequality \eqref{2.4} yields
\bea
\int\limits_\bbD |\varphi|^2  |f|^2 dA  
\geq \int\limits_G  |\varphi|^2  |f|^2 dA  
\geq S^2 \int\limits_G    |f|^2 dA 
> \frac{S^2}{M} \Vert f \Vert_2^2 =: \beta \Vert f \Vert_2^2 .  \label{2.6n}
\eea
Here, the constant $\beta > 0$  must be smaller than one, since the assumption 
$\Vert \varphi \Vert _\infty \leq 1$ implies
$
\Vert  f \Vert_2^2 \geq  \int_\bbD |\varphi|^2  |f|^2 dA .
$
Denoting by $I$ the identity operator on $A^2$, \ef{2.6a},  \ef{2.6n} yield
\bea
& & \Vert ( I - T_\varphi ) f\Vert_2^2
\leq \Vert M_{1 - \varphi} f \Vert_2^2
\leq \int\limits_\bbD |1 - \varphi|^2 |f|^2 dA 
\rowleq 
\int\limits_\bbD (1 - |\varphi|^2)  |f|^2 dA 
\leq (1 - \beta )  \int\limits_\bbD   |f|^2 dA ,   \label{2.6w}
\eea
hence, the invertibility of $T_\varphi$ follows from $\Vert I - T_\varphi \Vert \leq 1- \beta <1$.

As for the upper estimate for the operator norm, we obtain from the Neumann series 
\beas
\Vert T_\varphi^{-1} \Vert  \leq  \sum_{n=0}^\infty \Vert I-T_\varphi \Vert^n  
\leq  \sum_{n=0}^\infty (1- \beta)^n  = \frac1\beta =    \frac{M}{S^2}. \ \ \Box
\eeas

\BEL
\label{lem2.1a}
Let $\varphi : \bbD \to \bbC$ be a bounded, measurable function such that
\bea
{\rm Re} \, \varphi(z) \geq  | {\rm Im}\, \varphi(z) |^2
\label{2.5}
\eea
for all $z \in \bbD$. Then, for all sufficiently small constants  $R > 0$, the function
$\psi 
= R \varphi (z)$ satisfies
\bea
{\rm Re} \, \psi (z) \geq |\psi(z)|^2 . 
\label{2.5f}
\eea
\ENL

Proof. 
We first note that  that if the symbol $\varphi$ satisfies \ef{2.5}
and we define $\check\psi (z) = S \varphi(z)$ for some constant  $0 < S \leq 1$, then $\check \psi$ satisfies \ef{2.5},
too:
\bea
{\rm Re} \, \check \psi(z) = S \,{\rm Re} \, \varphi(z) \geq S
 \big( {\rm Im}\, \varphi(z)\big)^2 = S^{-1}  
 \big( {\rm Im}\,\check \psi(z)\big)^2 \geq  \big( {\rm Im}\, 
\check \psi(z)\big)^2 .  \label{2.5y}
\eea
We choose $S \in (0,1]$ so small that $\Vert \check \psi \Vert_\infty \leq 1$. Then, 
with the help of \ef{2.5y} we obtain 
\bea
& &  {\rm Re}\, \check \psi(z) =  
\frac12  {\rm Re}\, \check \psi(z)   +  \frac12{\rm Re}\, \check \psi(z) 
\geq 
\frac12  {\rm Re}\, \check \psi(z)   + 
\frac12  \big( {\rm Im}\, \check \psi(z)\big)^2 
\rowgeq 
\frac12  \big({\rm Re}\,\check \psi(z) \big)^2   + 
\frac12  \big( {\rm Im}\, \check \psi(z)\big)^2 
= \frac12  |\check \psi(z)|^2 .  \label{2.6}
\eea
Finally, if $R \in (0,1)$ is such that  $R \leq S /2$, we define $\psi = R \varphi$ 
so that there also holds $\psi = RS^{-1} \check \psi $. We obtain, by \ef{2.6} and the choice 
of $R$,
\beas
& & {\rm Re}\, \psi (z) = \frac{R}S  {\rm Re}\, \check \psi(z) \geq
\frac{R}{2S}  |\check \psi(z)|^2 = \frac{R}{2S} \Big( \frac{S}R \Big)^2
|\psi(z)|^2 
\roweq
\frac{S}{2} \cdot \frac1R |\psi(z)|^2 \geq |\psi(z)|^2  .
\ \ \Box
\eeas

\BER \label{rem2.2a}
We will need the observation that any number smaller than or equal to $ 2^{-1}
\min(1, \Vert \varphi \Vert_\infty^{-1})$ can be chosen for the number $R$, see the above proof. 
\ENR

\section{Main theorem and examples} \label{sec3}

Let us present the following generalization of Theorem 3.2. of \cite{ZZ16}.

\BET
\label{th2.2}
Let $\varphi : \bbD \to \bbC$ be a bounded, measurable function such that
\bea
{\rm Re} \, \varphi(z) \geq  |{\rm Im}\, \varphi(z)|^2 \ \ \ \ \forall \, z \in \bbD.
\label{2.9}
\eea
Then, the Toeplitz operator $T_\varphi$ is invertible on $A^2$, if
and only if the Berezin transform $B(|\varphi|)$ of the non-negative
symbol $|\varphi|$ is invertible in $L^\infty(\bbD)$. 
\ENT

Proof.  Theorem 3.1. of \cite{ZZ16} shows that the invertibility of the 
Toeplitz operator implies the invertibility of the Berezin transform.

Conversely, assume that the Berezin transform $B(|\varphi|)$ is invertible
in $L^\infty(\bbD)$. Note that for any constant $R>0$, the invertibility of the operator $T_\varphi$ 
is equivalent to the invertibility of $T_{R\varphi}$. Accordingly, we choose $R \leq 
\Vert \varphi \Vert_\infty^{-1}$ so small 
that \ef{2.5f} holds for $\psi = R \varphi$. Then, we  redefine $\varphi := \psi$ and note that the 
Berezin transform $B(|\varphi|)$ for the new $|\varphi|$ is still invertible
and there holds $\Vert \varphi \Vert \leq 1$. 

We proceed by using the proof of Theorem 3.2. of \cite{ZZ16} to the non-negative symbol 
$| \varphi|$: the proof shows that the function $|\varphi|$ satisfies $(iii)$ of Theorem \ref{lem1.1}.
By the statement $(ii)$ of  Theorem \ref{lem1.1} for the
function $|\varphi|$, we thus find a constant $ \eta > 0$ such that
\bea
\int\limits_\bbD |\varphi |^2  |f|^2 dA \geq \eta \int\limits_\bbD  |f|^2 dA   \label{2.9e}
\eea
for all $f \in A^2$. Now, if we choose $S > 0$ such that $S < \eta$ and define the set
$G = \{ z \, : \, |\varphi(z)| >  S\}$, 
we obtain by \ef{2.9e}
\beas
& & \int\limits_G |\varphi |^2  |f|^2 dA  = \int\limits_\bbD |\varphi |^2  |f|^2 dA
- \int\limits_{\bbD \setminus G}  |\varphi |^2  |f|^2 dA
\rowgeq 
\eta \int\limits_\bbD  |f|^2 dA  - \int\limits_{\bbD \setminus G}  S  |f|^2 dA 
\geq 
( \eta -S) \int\limits_\bbD  |f|^2 dA .
\eeas
Thus,  condition \ef{2.4} in Lemma \ref{lem2.1} of the  present paper holds with any 
constant $M > 1/(\eta -S)$. Condition \ef{2.5d} is also satisfied, by the scaling made in the 
beginning of the proof,  hence,
Lemma \ref{lem2.1} implies the invertibility of $T_\varphi$. \ \ $\Box$

\bigskip

We remark that condition \ef{2.9} in Theorem \ref{th2.2} can be replaced by the seemingly weaker
\bea
{\rm Re} \, \varphi(z) \geq \delta |{\rm Im}\, \varphi(z)|^2
\label{2.99}
\eea
where $\delta \in (0,1]$ is an arbitrary constant. Namely, if $\varphi$ satisfies
\ef{2.99}, then the symbol $\delta \varphi$ satisfies \ef{2.9}, and the invertibility properties are the same
for the symbol $\varphi$ and $\delta \varphi$. 

\bigskip

{\bf Example.} 
The values of the harmonic polynomial  $P(z) = 2+ \frac12 z + \frac12 \bar z^2  $
are in the parabolic region \ef{2.9}: 
the real part 
of $P(z) = P(x +iy)$ equals
\bea
2 + \frac12 x + \frac12 (x^2- y^2) 
\eea
On the circle  $y^2 = 1 - x^2$ this equals $\frac32 + \frac{x}2 + x^2$, which has minimum 
$1 \frac7{16} > \frac54$ at the  point $x = -1/4$.

The square of the imaginary part $\frac12 y - xy = (\frac12 - x)y$ 
boils down on the circle $y^2 = 1 - x^2$ to  $(1- x^2) (x-\frac12)^2$. This function
has its maximum at the point $x_0 = \frac18 - \frac14 \sqrt{\frac{33}4 }$, which satisfies
$-\frac58 < x_ 0 < -\frac48 = -\frac12$, and we obtain the estimate
\beas
& & (1- x_0^2) \Big( x_0-\frac12 \Big)^2 < \Big( 1 - \Big( \frac12\Big)^2 \Big) \Big( -\frac58- \frac12 \Big)^2
< \frac34 \Big(\frac98 \Big)^2 < 1 . 
\eeas
By the previous estimate, the inequality \ef{2.9} thus holds on the circle $\partial \bbD$ and consequently 
also on $\bbD$. 

Obviously,
\beas
|P(z)| \geq 2 - \frac12 |z|  - \frac12 |z^2 |  \geq 1 \ \ \ \ \forall\,z \in \bbD 
\eeas
so that the Berezin transform of $|P|$ is invertible. Theorem \ref{th2.2}
thus implies that the operator $T_P$ is invertible.

\bigskip

We next consider a general class of examples, which are harmonic polynomials. The following
observation will be needed.

\BER If the values of a mapping  $\varphi: \bbD \to \bbC$ satisfy
\bea
| 1 - \varphi(z)| < 1    \label{(1)}
\eea
for all $z \in \bbD$, then also \ef{2.99} holds with $\delta = 1/2$. Thus, \ef{2.9} of 
Theorem \ref{th2.2} is satisfied by the symbol $\frac12 \varphi $. 

To see this, if \ef{(1)} holds, we 
have $\varphi(z) - 1  = re^{i \theta}$ with  $0 \leq r < 1$, for all $z\in \bbD$.
We get
\bea
{\rm Re } \varphi(z) =   1 + {\rm Re }\big( \varphi(z) - 1 \big)
= 1+  r  \cos \theta  \label{(12)}
\eea
and 
\beas
& & | {\rm Im } \varphi(z)|^2  =   | {\rm Im }\big( \varphi(z) - 1 \big)|^2 
=  r^2  \sin^2 \theta = r^2 ( 1 - \cos^2 \theta) 
\roweq
r^2  ( 1 +  \cos \theta) ( 1-  \cos \theta) \leq
2 r^2   ( 1 +  \cos \theta) 
= 2 r (r +  r\cos \theta) \leq 2 r (1 +  r\cos \theta). 
\eeas
By \ef{(12)}, this is not larger than $2{\rm Re } \varphi(z)$, which implies \ef{2.99}. 
\ENR

In the next theorem we consider symbols, which are harmonic polynomials
\begin{equation}
P(z) = p_0 + \sum_{m \in M} p_m z^m + \sum_{n \in N} q_n \bar z^n ,  \label{(2)}
\end{equation}
where $p_0 \in \bbR$ with $p_0 \geq 1$ and $M, N \subset \bbN_1 = \{1,2,3,\ldots\}$ are finite sets 
and $p_m, q_n \in \bbC$ are {\it non-zero} coefficients for all $m\in M$, $n \in N$ such that
\begin{equation}
\sum_{m \in M} |p_m|  + \sum_{n \in N} |q_n |  = 1 .  \label{(3)}
\end{equation}

\BET
\label{th3.3}
Let the harmonic polynomial $P$ be as in \ef{(2)}--\ef{(3)}.

\noindent
$(i)$ If $p_0 > 1$, then the Toeplitz operator $T_P$ is invertible. 

\noindent $(ii)$ If $p_0 =1$, then $T_P$ is not invertible, if and only if the following 
condition  holds: there exist $\lambda \in [0, 2 \pi]$ and, for all $m\in M$ and $n \in N$,
integers  $k_m$ and $\ell_n$ such that
\begin{equation}
m \lambda + {\rm arg}\, p_m = \pi + k_m 2 \pi \ \ \ \mbox{and} \ \ \ 
 - n \lambda  + {\rm arg}\, q_n = \pi + \ell_n  2 \pi . \label{(6)}
\end{equation}
\ENT 

Here, arg\,$z$ denotes the principal value of the argument of $z$. 
Before the proof, let us consider some illuminating 

\bigskip

{\bf Examples.} 
$(a)$ If  $R(z) = 1+ \frac23 z^3 + \frac13 \bar z^3  $, then $T_R$ is not invertible.  
$(b)$ If $Q(z) = 1+ \frac23 z^2 + \frac13 \bar z^3  $, then $T_Q$ is invertible.

Let us prove these claims using item $(ii)$. 
As for $R$, in formula \ef{(2)} we have $M=N = \{3\}$, and to see \ef{(6)} we choose $\lambda = \pi/3$, $k_3 = 0$, $\ell_3 = -1$. Hence, with $m=3 =n$ we get
\beas 
m \lambda + {\rm arg}\, p_m = 3 \lambda  = \pi = \pi +k_3 2 \pi  \ , \ \ - n \lambda +
{\rm arg}\, q_n = - 3 \lambda = - \pi = \pi + \ell_3 2 \pi,
\eeas
so that \ef{(6)} holds and thus $T_R$ is not invertible. 

Concerning  $Q$, in formula \ef{(2)} we have $M=\{2\}$, $N=\{3\}$. Suppose $\lambda \in [0, 2\pi]$
and $k_2 , \ell_3 \in \bbZ = \{0, \pm1, \pm 2, \ldots\}$ are such that \ef{(6)} holds. Then, ${\rm arg}\, p_2
= 0 = {\rm arg}\, q_3$ and  thus
\beas
\lambda = \frac{\pi}{2} + k_2 \pi = \pi \Big( \frac12 + k_2 \Big)
\ \ \mbox{and} \ \ \lambda = - \frac{\pi}{3} - \ell_3 \frac{2 \pi}3
= \pi \Big( - \frac13 - \frac{2 \ell_3}3 \Big) 
\eeas
so that equating both expressions for $\lambda$, we get 
\bea
\frac12 + k_2  = - \frac13 - \frac{2 \ell_3}3 , \label{(7)}
\eea
But there does not exist integers  $k_2 , \ell_3$  such that \ef{(7)} holds:
on the right, the number is always integer $+ 0, 1/3$ or $2/3$, so it cannot be integer + 1/2.
Hence, \ef{(6)} cannot be satisfied. As a conclusion, $T_Q$ is invertible.

\bigskip

Proof of Theorem \ref{th3.3}. To prove the claim $(i)$, we 
write $p_0 = 1 + q_0 $, where $q_0 > 0$, and consider first
the harmonic symbol 
\bea
\varphi(z) = 1 + \sum_{m \in M} p_m z^m + \sum_{n \in N} q_n \bar z^n ,  \label{(4)}
\eea
i.e.  $P = \varphi + q_0$. 
Then, we have, by the triangle inequality and \ef{(3)}, 
\bea
& & | 1- \varphi(z)| = \Big| \sum_{m \in M} p_m z^m + \sum_{n \in N} q_n z^n \Big|
\leq \sum_{m \in M} |p_m| |z|^m + \sum_{n \in N} |q_n| |z|^n
\row<
\sum_{m \in M} |p_m|  + \sum_{n \in N} |q_n| = 1  \label{(5)}
\eea
so that \ef{(1)} and thus \ef{2.99} hold for $\varphi$. Since
${\rm Re}\, P(z) = {\rm Re} \,  \varphi (z) + q_0 > 
{\rm Re} \,  \varphi (z)$ and ${\rm Im} \, P(z) =  {\rm Im} \, 
\varphi (z)$, \ef{2.99} holds also for $P$.

To show that the  Berezin transform of $|P|$  is invertible on $\bbD$, we again
use  \ef{(3)} and the triangle inequality to obtain
\beas
|P(z)| &\geq &  p_0 -  \Big| \sum_{m \in M} p_m z^m + \sum_{n \in N} q_n z^n \Big|
\geq p_0 - \bigg(\sum_{m \in M} |p_m|  +  \sum_{n \in N} |q_n|\bigg)
\roweq
p_0 -1 = q_0 > 0 
\eeas
for all $z \in \bbD$. The constant function $z \mapsto q_0$ is harmonic, hence $B(q_0) =
q_0$ so that the above estimate proves $B(|P|)(z) \geq q_0$ for $z \in \bbD$. Thus $B(|P|)$ is invertible,
and the claim $(i)$ follows from Theorem \ref{th2.2}. 

We turn to claim  $(ii)$. Let $P$ be as in \ef{(1)}--\ef{(2)} with $p_0=1$. 
Note that now $P$ coincides with $\varphi$ in the proof of $(i)$, and hence, by
\ef{(5)}, $P$ satisfies \ef{(1)} and thus also \ef{2.99}. Thus, Theorem \ref{th2.2}
applies and we need to consider the invertibility of the Berezin transform
$B(| P|)$ as a function in $L^\infty(\bbD) $.

Let us first prove the "if"-part: we assume that one can find $\lambda$, $k_m$, $\ell_n$ 
such that \ef{(6)} holds and show that $B(|P|)$ is not invertible. Indeed, 
using \ef{(6)} and \ef{(3)}  we get
\bea
P(e^{i \lambda} ) &=& 1 + 
\sum_{m \in M} p_m e^{i m \lambda} + \sum_{n \in N} q_n e^{- i n \lambda} 
\roweq 
1 +  \sum_{m \in M} |p_m| e^{i ({\rm arg} \,p_m  +  m \lambda) } + \sum_{n \in N} |q_n| e^{i
( {\rm arg} \, q_n +  n \lambda) }
\roweq 
1 +  \sum_{m \in M} |p_m| e^{i \pi  } + \sum_{n \in N} |q_n| e^{i \pi} 
= 1 -  \sum_{m \in M} |p_m| - \sum_{n \in N} |q_n| = 0 . \label{(8)}
\eea
Thus, also $|P|$ has a zero on the boundary $\partial \bbD$, and since the values of 
$|P|$ and $B(|P|)$ coincide on $\partial \bbD$, we conclude that 
the Berezin transform of $|P|$ is not invertible
in $L^\infty$. Theorem \ref{th2.2} shows that $T_P$ is not invertible. 

To prove the ``only if"-statement we assume the numbers  $\lambda$, $k_m$, $\ell_n$ 
cannot be found for every $m,n$ so as to satisfy \ef{(6)}. Thus, if  $\lambda \in [0, 2 \pi]$ is 
arbitrary, there is $ m_\lambda \in M$ such that  the first identity in 
\ef{(6)} does not hold for any integer $k$ (or, the second identity in \ef{(6)} fails, which case
is treated in the same way.) 

We thus have 
$$
 m_\lambda \lambda  + {\rm arg} \, p_{m_\lambda}  \not= \pi + k 2 \pi 
$$
for all $k \in \bbZ$, which means that 
$$
e^{i ({\rm arg} \,p_{ m_\lambda}  +   m_\lambda \lambda) } \not= -1 .  
$$
Every point on $\partial \bbD$, except -1, has real part bigger than -1, so that we can write
\bea
{\rm Re}\, e^{i ({\rm arg} \,p_{m_\lambda}  +   m_\lambda \lambda) } = -1 + \delta(\lambda) \label{(9)}
\eea
for some $\delta(\lambda) > 0$. The compactness of the interval $[0, 2\pi]$ implies that there is a number $\delta' > 0$, independent of $\lambda$,  such that
\bea
{\rm Re}\, e^{i ({\rm arg} \,p_{m_\lambda}  +  m_\lambda \lambda) } \geq  -1 + \delta'. \label{(10)}
\eea
We postpone the proof of this observation to Lemma \ref{lemY}, below.

Given $z\in\bbD$ with $|z| \geq 1/2$ we write $\lambda = {\rm arg}\,z$ and obtain by \ef{(10)}, \ef{(3)}
and a repeated use of the triangle inequality
\bea
{\rm Re}\,  P(z) &=& 1 + {\rm Re}\,( p_{m_\lambda} z^{m_\lambda} )
+{\rm Re}\, \sum_{m \in M\setminus \{ m_\lambda\} } p_m z^m + {\rm Re}\,\sum_{n\in N} q_n \bar z^n
\roweq
1 + {\rm Re}\, e^{i ({\rm arg} \,p_{m_\lambda}  +  m_\lambda \lambda) }  |p_{m_\lambda}z^{m_\lambda}| 
+{\rm Re}\, \sum_{m \in M\setminus \{ m_\lambda\}} p_m z^m + {\rm Re}\,\sum_{n\in N} q_n \bar z^n
\rowgeq
1 + \delta' |p_{m_\lambda}z^{m_\lambda}|  - |p_{m_\lambda}z|^{m_\lambda} 
-  \sum_{m \in M\setminus \{ m_\lambda\}} |p_m z^m| - \sum_{n\in N} | q_n \bar z^n |
\rowgeq
1 + \delta' |p_{m_\lambda}z^{m_\lambda} | -  \sum_{m \in M} |p_m | - \sum_{n\in N} | q_n  |
\geq \delta' |p_{m_\lambda}|2^{- m_\lambda} 
\rowgeq 
\delta' 2^{-\mu} \min\limits_{m \in M } \{ |p_m|\}, \label{(11)}
\eea
where $\mu = \max_{m \in M} \{m \}$. The number on the right of \ef{(11)}  is positive and does not 
depend on $z$. For $|z| < 1/2$ we obtain
from \ef{(2)} and \ef{(3)}
\beas
|P(z)| \geq 1 - \sum_{m \in M}  |p_m| 2^{-m} -  \sum_{n \in N} |q_n| 2^{-n}
\geq 1 - \frac12 \bigg( \sum_{m \in M}  |p_m| +  \sum_{n \in N} |q_n|\bigg)
\geq \frac12  . 
\eeas
This and \ef{(11)} show that $|P(z)|$ is bounded from below by a positive constant $\delta > 0$. 
This implies that  $B(|P|)(z) \geq \delta$ for all $z \in \bbD$, hence
$B(|P|)$ is invertible as an element of $L^\infty(\bbD)$. By Theorem 3.1, also $T_P$ is invertible. 
\ \ $\Box$

\bigskip

Let us finally complete the missing step in the proof of Theorem \ref{th3.3}. We keep the same 
notation as above.

\BEL
\label{lemY}
There exists a number $\delta' > 0$  such that
the inequality \eqref{(10)} is satisfied for all $\lambda \in [0,2\pi]$. 
\ENL

Proof. For every $\lambda \in  [0, 2\pi]$, the identity \ef{(9)} implies that 
$$
{\rm Re}\, e^{i ({\rm arg} \,p_{m_\lambda}  +   m_\lambda \nu) } \geq -1 + \delta(\lambda)/2 
$$
for $\nu$ belonging to a small interval $I_\lambda := \big(\lambda - \epsilon(\lambda) , \lambda + 
\epsilon(\lambda)\big)$. The intervals $I_\lambda$ form an open covering of $[0,2\pi]$, so that by 
compactness, we can
pick up a finite subcovering, and then define $\delta'$ to be the smallest of the corresponding 
numbers $\delta (\lambda)$. Every $\lambda \in [0, 2 \pi]$ belongs to some of these finitely many
intervals $I_{\lambda_j}$, $j= 1, \ldots, J$, and we finally
redefine $m_\lambda \to m_{\lambda_j}$ and $p_{m_\lambda} \to p_{m_{\lambda_j}}$ for 
$\lambda \in I_{\lambda_j}$; a number $\lambda$ may belong to more than one $I_{\lambda_j}$,
but we just choose one of these for every $\lambda$. We thus get \ef{(10)} for all $\lambda \in  [0, 2\pi]$.
\ \ $\Box$

\bigskip

\section{Sufficient conditions for the invertibility} \label{sec7}

In this section, Theorem \ref{th3.1}, we show that  the  assumption \ef{2.9} can still be weakened by assuming 
it only for $z$ bounded away from zero  (see $(i)$, $(ii)$), although one then needs another condition 
$(iii)$ 
to control the range of $\varphi$ near the origin.  Note that the result allows $\varphi(z)$ to have a negative
real part for some $z\in \bbD$. 

We first recall a known norm bound for the multiplication operator. 
Denote $e_n = \sqrt{2(n+1)} z^n$ for all $n \in \bbN$ so that these functions form an orthonormal basis of $A^2$. 
Let  us denote, for all $n \in \bbN$ by $A_n$ the $n$-codimensional 
closed subspace of $A^2$  spanned by monomials of degree at least $n$.
Now, if $\omega \subset B(0,r) = \{ z \, : \, |z| \leq r \}$ for some $0 < r < 1$
 and $\chi_\omega$ is its characteristic (indicator)
function  and $M_\omega := M_{\chi_\omega}$, then the operator norm 
of $M_\omega$ satisfies
\bea
\Vert M_\omega \big|_{A_n} \Vert \leq \frac{ r^{n+1}}{(1-r)^{1/2}}. \label{lem3.2}
\eea
To see this, we have for all $m \geq n$
\bea
\int\limits_\bbD |M_\omega z^m|^2 dA \leq \int\limits_0^r s^{2m +1} ds 
= \frac{1}{2m +2} r^{2m +2}
\eea
hence, for $f = \sum_{m \geq n} a_m e_m \in A_n$ we get by the Cauchy-Schwartz inequality
\bea
& & \Big\Vert M_\omega \sum_{m=n}^\infty a_m e_m \Big\Vert_ 2
\leq \sum_{m=n}^\infty  |a_m | \, \Vert M_\omega e_m \Vert_2
\leq \sum_{m=n}^\infty  |a_m | r^{m+1} 
\rowleq
\Vert f \Vert_ 2 \Big( \sum_{m=n}^\infty  r^{2(m+1)} \Big)^{1/2}
\leq
\frac{ r^{n+1}}{(1-r)^{1/2}}
\Vert f \Vert_ 2  .
\eea

For a set $A \subset \bbC$ we denote by
$|A| $ its area.

\BET
\label{th3.1}
Let $\varphi \in L^\infty$ and let 
the number $\rho \in (0, \rho_0)$ be arbitrary, where $\rho_0 $ is the constant $ \min(1, \Vert \varphi \Vert_\infty^{-1})/32$.  If

\medskip

\noindent $(i)$  
$ 
{\rm Re} \, \varphi(z) \geq  |{\rm Im}\, \varphi(z)|^2 
$
for all $ z\in \bbD$ with $ |\varphi(z)| \geq \rho$, 

\medskip

\noindent $(ii)$ $|\varphi (z)| > \rho $ for all $z$ with $|z| \geq \rho$,  and

\medskip

\noindent $(iii)$  the set  $\Lambda := \{ z \in \bbD\, : \, |\varphi(z)| \leq \rho \}$
satisfies  $|\Lambda | \leq \rho^6$, 

\medskip

\noindent
then the Toeplitz operator  $T_\varphi$ is invertible on $A^2$.
\ENT

Below 
we also consider symbols $\varphi$ for which the power $\rho^6$ in $(iii)$ can be improved into $\rho^4$.  

\bigskip

Proof. 
Given $\varphi$,  we first consider 
an arbitrary  $\rho \in (0,1) $, assume that  $(i)$--$(iii)$ hold for $\rho$ and finally show that 
$T_\varphi$ is  invertible, if $\rho \in (0,\rho_0)$ for a small enough $\rho_0$. 

Let us  denote by  $\chi : \bbD \to 
[0,1]$ the  characteristic function of the subdomain $\Omega= \bbD \setminus \Lambda = 
\{ z \in \bbD \, : \, |\varphi(z)| > \rho\}$. Then, by $(i)$ and the definition of the set $\Lambda$,
the function $\chi \varphi$ satisfies condition \ef{2.5} of Lemma \ref{lem2.1}. By Remark
\ref{rem2.2a} and the bound $\Vert \chi \varphi \Vert_\infty \leq \Vert \varphi \Vert_\infty$,
\ef{2.5f} holds for the function $R \chi \varphi$, when we choose $R:= 2^{-1} \min(1, 
\Vert  \varphi \Vert_\infty^{-1}) \leq 1/2 $. (Note that $R$ does not depend on $\rho$, although $\Omega$ 
and $\chi$ do.)  In other words, the inequality 
\bea
{\rm Re}\, \psi(z) \geq |\psi(z)|^2 \label{2.6af}
\eea
holds for the function $\psi = R \varphi$, for all $z \in \Omega$. (Note that $\Vert 
\psi \Vert_\infty \leq 1$.)  As in \ef{2.6a}, we then conclude that 
\bea
|1 - \psi(z) |^2 \leq 1 - |\psi(z)|^2 \ \ \ \forall \, z \in \Omega. \label{2.6aa}
\eea

We denote $S= 1/R$ and $\sigma = R \rho \in (0,1)$. For an arbitrary $f = \sum_{n=0}^\infty f_n z^n \in A^2$ 
we obtain 
\bea
(SI- T_\varphi) f = (SP- T_\varphi)f =
(SP- T_\varphi)(\chi f + (1- \chi)f )  \label{3.10}
\eea
Then, we take into account  $|\psi(z)| = R |\varphi(z) | \geq R \rho 
= \sigma$ for all  $z \in \Omega = \bbD \setminus  \Lambda$ and thus can  estimate
\bea
& & \Vert (S P  - T_{ \varphi} )\chi f\Vert_2^2
\leq \Vert (S I - M_{ \varphi} )\chi f\Vert_2^2
=   \Vert (S I - S M_{ \psi} )\chi f\Vert_2^2
\roweq
S^2 \int\limits_\Omega |1 - \psi|^2 |f|^2 dA 
\leq S^2  \int\limits_\Omega (1 - |\psi|^2)  |f|^2 dA 
\leq 
S^2  (1- \sigma^2) \int\limits_\Omega  |f|^2 dA  .
  \label{3.12}
\eea
where also \ef{2.6aa} was used for $\psi$.

We write $g(z) = f_0 + f_1 z$ and $h(z) = f (z) - g  (z)
= \sum_{n \geq 2 }f_nz^n $  and proceed with
\bea
& & \Vert (S P - T_{ \varphi} )(1-\chi) f\Vert_2
\rowleq 
\Vert (S P - T_{ \varphi} )(1-\chi) g\Vert_2
+ \Vert (SP - T_{ \varphi} )(1-\chi) h\Vert_2 .
\eea
By $(ii)$, the set $\Lambda$ is contained in $B(0,\rho) \subset \bbC$, hence \ef{lem3.2} and the 
definition of $\chi$ give us
\bea
& & \Vert (S P - T_{ \varphi} )(1-\chi) h\Vert_2^2
\leq \Vert (S  - SM_{ \psi} )(1-\chi) h\Vert_2^2
=  S^2 \int\limits_{\Lambda} |1 - \psi|^2 |h|^2 dA 
\rowleq 
S^2 \int\limits_{B(0, \rho)}  |h|^2 dA 
= S^2 \Vert M_{B(0,\rho)} h \Vert_2^2 \leq  
\frac{ S^2 \rho^6}{1- \rho } \Vert h \Vert_2^2
\leq  \frac{ S^2\rho^6}{1- \rho }  \Vert f \Vert_2^2 .  \label{3.14}
\eea
Moreover,  by $(iii)$, 
\bea
& &  
\Vert (S P - T_{ \varphi} )(1-\chi) g\Vert_2^2 \leq
\Vert (S I - M_{ \varphi} )(1-\chi) g\Vert_2^2
\rowleq \int\limits_{\Lambda} |S - \varphi(z)|^2 (|f_0|^2 + |f_1|^2) dA(z)
\rowleq 
4(S + \rho)^2 |\Lambda| \, \Vert f \Vert_2^2 \leq 
4(S + \rho)^2 \rho^6 \Vert f \Vert_2^2 
 \label{3.16}
\eea

Combining \ef{3.10}--\ef{3.16}  and taking into account  
$\sqrt{1 - \sigma^2} \leq 1 -  \sigma^2 / 2 $ we get
\bea
& & \Vert ( SI -T_\varphi )f \Vert \leq \big( S \sqrt{1- \sigma^2} + 
S  \rho^3 (1 -\rho)^{-1/2} + 2 ( S + \rho)\rho^3 \big) \Vert f \Vert_2
\rowleq
S \Big( 1 - \frac{R^2 \rho^2}{2}   +  \rho^3\Big(\frac1{\sqrt{1-\rho}} +
2+ 2\rho R  \Big) \Big) \Vert f \Vert_2  .
\label{3.17}
\eea
Note that due to the choice of $\rho_0$ and $R= 2^{-1} \min(1, \Vert  \varphi \Vert_\infty^{-1})$,
we have
$ 
\rho < R^2 / 8 , 
$ 
hence,
\beas
\rho^3 \Big( \frac1{\sqrt{1-\rho}} + 2+ 2\rho R  \Big)
\leq \rho^2 \frac{R^2}2 \frac14\Big( \frac1{\sqrt{7/8}} + 2+ \frac14 \Big)
= \rho^2 \frac{ R^2 \sigma_1}2 
\eeas
with a constant  $\sigma_1  \in(0,1) $. One obtains for \ef{3.17} the bound
\bea
S \Big( 1 - \frac{ R^2 \rho^2 }{2} ( 1 - \sigma_1  ) \Big) \Vert f \Vert_2 ,
\eea
which  yields  $\Vert SI - T_\varphi \Vert \leq \sigma_2 S $ for another positive constant 
$\sigma_2 < 1$. Hence, 
\bea
\Vert I - T_{R \varphi} \Vert = R  \Vert SI - T_\varphi \Vert \leq 
RS \sigma_2 = \sigma_2 < 1, 
\eea
and we obtain from the Neumann series
that $T_{R \varphi}$ and thus also $T_\varphi$ are invertible. 
\ \ \ 
$\Box$ 

\bigskip

The proof of Theorem \ref{th3.1} can be modified to get other similar statements.
In particular, the stronger requirement
$ 
{\rm Re} \, \varphi(z) \geq  |{\rm Im}\, \varphi(z)| ,
$
would allow us to replace condition $(iii)$ by the weaker assumption that
the set  $\Lambda := \{ z \in \bbD\, : \, |\varphi(z)| \leq \rho \}$
satisfies  $|\Lambda | \leq \rho^4$. We leave the details for the interested reader.

The following sufficient condition for the invertibility of $T_\varphi$ follows from 
Proposition \ref{th3.1}. We remark that condition \ef{2.99} could be used in the place of the first inequality of \ef{3.30}.

\BEC \label{cor3.5} 
Assume $\varphi:\bbD \to \bbC $ is  bounded, measurable and that for some $\delta  
\in (0,1]$  there holds
\bea
{\rm Re} \, \varphi(z) \geq  |{\rm Im}\, \varphi(z)|^2  \ \ 
\mbox{and} \ \   \delta \leq  |\varphi(z) | \leq \Vert \varphi \Vert_\infty 
\label{3.30}
\eea
for all $z \in \bbD$.
Then, the Toeplitz operator $T_\varphi$ is invertible on $A^2$. 
\ENC

Proof. 
If $\varphi$ is as in \ef{3.30}, it also satisfies $(i)$  
for any $\rho$ and $(ii)$ for all $\rho$ smaller than $\delta$ of \ef{3.30}. Moreover, if $\rho < 
\delta$, then the set $\Lambda$   is empty and $(iii)$ automatically holds. 
Thus, Proposition \ref{th3.1} implies the invertibility of $T_\varphi$.  \ \ $\Box$

\bigskip

We finally remark that in the case of a harmonic symbol $\varphi$, Corollary \ref{cor3.5} yields a 
partial positive answer to Douglas' question for Bergman spaces, since the latter condition in \ef{3.30} 
coincides with assuming the invertibility of the Berezin transform of $\varphi$.

\section{Results concerning harmonic symbols, Berezin transform and its iterations} \label{sec6}

Next we consider Toeplitz operators with harmonic symbols.  
Let $\cP : L^\infty(\partial \bbD) \to L^\infty(\bbD)$  denote the Poisson extension of a bounded
function on the circle to a bounded harmonic function on $\overline{\bbD}$. We first present 
a reformulation of Theorem \ref{th2.2} in this setting. Note that the assumptions in the
next proposition are the same as in Douglas' theorem, except for condition \ef{h.2}. 

\BEP  \label{prop5.1} 
Assume $g: \partial \bbD \to \bbC$ is continuous and satisfies
the condition 
\bea
{\rm Re}\, g(e^{i \theta}  ) \geq  | {\rm Im} \, g(e^{i \theta}  )|^2 
\ \ \ \ \forall \, \theta \in [0, 2 \pi]. \label{h.2}
\eea
Then, the Toeplitz operator with symbol $\varphi =  \cP g$ is invertible on $A^2$, if and only if
$g(e^{i \theta}  ) \not= 0$ for every $\theta$. 

\ENP

Proof. Assume first $g(e^{i \theta}  ) \not= 0$ for every $\theta$. Since \eqref{h.2} holds, there cannot
exist a point $t \in [0, 2 \pi]$ such that ${\rm Re} \, g(e^{i t } ) =0$ but ${\rm Im} \, 
g(e^{i t}  ) \not= 0$. Since the case ${\rm Re} \, g(e^{i t } ) =  {\rm Im} \, g(e^{i t}  )= 0$ is now also 
excluded, we conclude  that ${\rm Re} \, g(e^{i t } ) >0 $ and thus,  by the continuity of $g$,  
${\rm Re} \, g(e^{i t } ) \geq c $ for some constant $c>0$. 
Then, the continuity of $\cP g$ in the closed disc $\overline \bbD$ and the positivity of the  Poisson kernel  on $
\bbD$ implies that   ${\rm Re} \, \varphi(z) =  {\rm Re} \,  \cP g(z) \geq \delta'$ for some constant 
$\delta' > 0$, for all $z \in \bbD$. Hence,
\beas
{\rm Re}\, \varphi(z)  \geq  \delta | {\rm Im} \, \varphi(z)|^2 \ \ \ \ \forall\, z \in \bbD,
\eeas
where $\delta = \delta' \Vert \varphi \Vert_\infty^{-2}$. By the remark around \ef{2.99},
the symbol $\delta \varphi$ satisfies condition \ef{2.9} of Theorem \ref{th2.2}.  Moreover, the Berezin
transform of $\delta \varphi$ coincides with $\delta \varphi$ and is thus invertible in $L^\infty(\bbD)$. 
Theorem \ref{th2.2} implies that $T_{\delta \varphi}$ and thus also $T_\varphi$ are invertible. 

If  $g(e^{i \theta}  ) = 0$ for some $\theta$, then $\varphi$ and thus $ B (| \varphi|)$
cannot be invertible in $L^\infty(\bbD)$. By Theorem 3.1. of \cite{ZZ16}, $T_\varphi$ is not
invertible. \ \ $\Box$ 

\bigskip

The following result yields many examples of invertible Toeplitz operators. 
We denote for all $n \in \bbN$ by  $B^n$ the $n$th iterate of the Berezin transform.

\BEP
Let $\varphi \in C(\overline \bbD)$ be such that $\big| B(\varphi) (z) \big| \geq \delta$ for some
constant $\delta > 0$ and all $z \in \bbD$.  

\noindent $(i)$ For every $n \in \bbD$, the Toeplitz operator $T_\psi$ with symbol $\psi = 
\cP\big (|\varphi|^n\big|_{\partial \bbD} \big)$
is invertible in $A^2$. 

\noindent $(ii)$ If $n \in \bbN$ is such that
\bea
{\rm Re}\, B^n \varphi(z)  \geq  \delta | {\rm Im} \, B^n\varphi(z)|^2 \ \ \ \ \forall\, z \in \bbD,
\label{6.4}
\eea
then the Toeplitz operator $T_\psi$ with symbol $\psi = B^n \varphi$ is invertible in $A^2$. 

\ENP

Proof. $(i)$ Since $\varphi \in C(\overline \bbD)$, the values of $\varphi$ and $B \varphi$  coincide on $
\partial \bbD$. The assumption of the proposition implies that $|\varphi|$ and thus $|\varphi|^n$ are 
non-zero on $\partial \bbD$. Denoting the restriction $g = |\varphi|^n\big|_{\partial \bbD}$, we find that $g$ 
satisfies the assumptions of Proposition \ref{prop5.1}, which implies the statement (note that \ef{h.2} holds, 
since $g$ is non-negative). 

$(ii)$ 
Let us fix  $n \in \bbD$ such that \ef{6.4} holds. Since $\varphi \in C(\overline \bbD)$, the 
values of $\varphi$ and every $B^k \varphi$ with  $k \in \bbN$ coincide on $\partial \bbD$, by an 
induction argument. 
Let us define the continuous function  $g = \varphi\big|_{\partial \bbD} = B^n \varphi\big|_{\partial \bbD}$  
on $\partial \bbD$. The continuity of the function $B^n \varphi$ on $\overline{\bbD}$ and \ef{6.4} imply that 
condition \ef{h.2} holds for $g$. Moreover, by the assumption $\big| B(\varphi) (z) \big| \geq \delta$, the function $g$ is  non-zero on $\bbD$. 
Since the integral kernel of the linear operator $B^n$ is positive, the same argument as in 
the proof of Proposition \ref{prop5.1} shows that there exists a constant $\delta  > 0$ 
such that $ {\rm Re}\, \varphi(z) \geq \delta  $ and thus also $   | B^n \varphi(z)| \geq \delta$
for all $z \in \bbD$. The invertibility of the operator $T_\psi$ now follows from Corollary 
\ref{cor3.5}. \ \ $\Box$

\bigskip

Finally, we consider the relation of the condition \eqref{2.5} for the symbol itself
and for its Berezin tranform. 

\BEP \label{prop6.3}
If the function  $\varphi \in L^\infty(\bbD)$ satisfies \eqref{2.5}, then $B( \varphi)$  
also satisfies the same condition.  
\ENP

Proof. Indeed, there holds 
\bea
B(\varphi)(z) = \int\limits_\bbD |k_z(w)|^2 {\rm Re} \varphi \, dA(w)
+ i \int\limits_\bbD |k_z(w)|^2 {\rm Im} \varphi \, dA(w).
\eea
so that using the Jensen inequality we obtain
\beas
& &  \Big(  {\rm Im}  \int\limits_\bbD |k_z(w)|^2\varphi \, dA(w)  \Big)^2
\leq 
 \Big( \int\limits_\bbD |k_z(w)|^2 \big| {\rm Im} \varphi \big| \, dA(w) \Big)^2 
\rowleq 
 \int\limits_\bbD |k_z(w)|^2 \big( {\rm Im} \varphi \big)^2 \, dA(w)
\leq 
\int\limits_\bbD |k_z(w)|^2  {\rm Re} \, \varphi \, dA(w)
\roweq {\rm Re} \, \int\limits_\bbD |k_z(w)|^2  \varphi \, dA(w) .  \hskip2cm \Box
\eeas

Proposition \ref{prop6.2} shows that the converse of implication  does not hold in general.

\BEP
\label{prop6.2}
There exist a (real valued, radial) symbol $\varphi$  such that the Berezin transform $B(\varphi)$ satisfies condition \ef{2.5} in the place of  $\varphi$ and is invertible, 
but $T_\varphi$ is not. Hence, condition \ef{2.5} does  not hold for the symbol
$\varphi$ itself. 
\ENP

Proof. In the proof of Corollary 3.5 of \cite{ZZ16} it is shown that  for 
$ P(z) = |z|^2 - \frac32 |z| + 1$ there holds 
\bea
\inf_{z \in \bbD} B(T_P)(z) > \frac{13}{28}  =: \lambda_2   \label{6.2}
\eea
where $\lambda_2$ is the lowest eigenvalue of $T_P$, and it  
corresponds to the eigenfunction  $f(z) = z^2$. Consequently,
\bea
Q(z) = P(z) - \frac{13}{28} = |z|^2 - \frac32 |z|
+ \frac{15}{28}    \label{6.3}
\eea
is a real valued, radially symmetric symbol with positive Berezin transform, but
$T_Q$ is not invertible, since it has 0 as an eigenvalue. 

In view of Theorem \ref{th2.2}, the symbol $\varphi$ cannot satisfy condition \ef{2.5}.
As an alternative proof, one can calculate that
the function $r \mapsto Q(r)$ has two real zeros in the interval $(0,1) \ni r$
so that $Q$ attains negative values on $\bbD$ and thus \ef{2.5} cannot hold. \ \ $\Box$

\bigskip


Finally, we show that if a symbol $\varphi$ satisfies condition \ef{2.9} and has some other general 
regularity properties, then the operator $T_\varphi$ is invertible, if and only if it is Fredholm.
This result follows from the next lemma.

\BEL
\label{lemX1}
Assume $\varphi \in C(\overline \bbD)$ satisfies condition \ef{2.9}. If there exists $\delta > 0$ 
such that $|B(\varphi)(z) | \geq \delta$ for all $z \in \partial \bbD$, then $B(\varphi)$ is invertible
as an element of $ L^\infty(\bbD)$. 
\ENL

Proof. Without loss of generality, we  assume that $\Vert \varphi \Vert_\infty \leq  1$, since this can be 
achieved by multiplying  $\varphi $ by a small enough number belonging to $(0,1)$, and it suffices to prove the 
lemma for this scaled symbol. 
By for example \cite{Zh}, Proposition 6.14, the Berezin transform $B(\varphi)$ is continuous
on $\overline \bbD$ so that, for some $\delta_1 > 0$, the condition $|B(\varphi)(z) | \geq \delta_1$ holds for all 
$z \in \bbD$ with $|z| \geq 1 - \delta_1$. Since $B(\varphi)$ is defined by an integral (cf. \ef{1.10})), 
there must exist a closed subset $\Omega \subset \bbD$ with positive area such that $|\varphi(z)| \geq \delta_2$
for some number  $\delta_2 > 0$ and all $z \in \Omega$. Using  \ef{2.9}, we thus obtain
\bea
& & {\rm Re} \, \varphi(w) \geq \frac12 {\rm Re} \, \varphi(w) + \frac12 {\rm Re} \, \varphi(w)
\geq \frac12 {\rm Re} \, \varphi(w) + \frac12 {\rm Im} \, \varphi(w)^2
\rowgeq
\frac12 {\rm Re} \, \varphi(w)^2  + \frac12 {\rm Im} \, \varphi(w)^2
= \frac12 |\varphi(w)|^2 \geq \frac12 \delta_2^2 \label{x1}
\eea
for $w \in \Omega$. 

We also have 
\bea
\inf \big\{ |k_z(w)| \, : \, |z| \leq 1 - \delta_1, \ w \in \Omega \big\} \geq \delta_3  \label{x2}
\eea
for some $\delta_3>0$, due to the compactness of the $z$- and $w$-sets here. Noting that the 
real part of $\varphi$ is by assumption \ef{2.9} non-negative everywhere in $\bbD$ and using
\ef{x1} and \ef{x2}, we get for all $|z| \leq 1 - \delta_1$
\beas
& & |B(\varphi)(z)| \geq {\rm Re} \, B(\varphi)(z) = \int\limits_\bbD |k_z(w)|^2 \,  {\rm Re} \, \varphi(w)  dA(w)
\rowgeq  \int\limits_\Omega |k_z(w)|^2 \, {\rm Re} \, \varphi(w) dA (w)
\geq \delta_4
\eeas
for some number $\delta_4>0$. This proves the claim, since it was observed above that 
$|B(\varphi)|$ is bounded from below in the set $\{ z \, : \, |z| \geq 1- \delta_1\}$. 
 \ \ $\Box$

\bigskip
Let us denote by $AQ(\bbD)$ the closed subalgebra of $L^\infty(\bbD)$ consisting of functions $\psi$ 
such that the Hankel operator $H_\psi = (I - P)M_\psi$ is compact.

\BET
\label{hX2}
Assume that $\varphi \in C(\overline \bbD) \cap AQ(\bbD)$ 
satisfies condition \ef{2.9}. Then, the Toeplitz operator $T_\varphi$ is invertible, if and only if 
it is a Fredholm operator.
\ENT

Proof. Let us assume that $T_\varphi:A^2 \to A^2$ is a Fredholm operator for a symbol 
$\varphi \in C(\overline \bbD) \cap AQ(\bbD)$. It 
follows from  Stroethoff-Zheng -92, Toeplitz and Hankel operators on Bergman spaces, TAMS 329 (1992) 773--794, Corollary 22, that 
we have $|B(\varphi) (z) | \geq \delta$ 
for $z \in \partial \bbD$ and some constant $\delta > 0$. We obtain from Lemma \ref{lemX1} that
$B(\varphi)$ and thus also $|B(\varphi)|$ are invertible in $L^\infty(\bbD)$. Theorem \ref{th2.2}
yields that $T_\varphi$ is invertible. 

The other implication of the theorem is trivial. \ \ $\Box$

\bigskip

{\it Acknowledgements.} The authors wish to thank Dragan Vukotic (Madrid) for some useful remarks on the 
manuscript, which helped to improve the formulation of Theorem \ref{th2.2}.

{\it Funding sources.} The second named author was partially supported by the Magnus Ehrnrooth Foundation and 
the V\"ais\"al\"a Foundation of the Finnish Academy of Sciences and Letters.


\begin{thebibliography}{99}


\bibitem{AC} P.~Ahern,  $\check{\mbox{Z}}$.~$\check{\mbox{C}}$u$\check{\mbox{c}}$kovi\'c, A theorem of Brown–Halmos type for Bergman space Toeplitz operators, J. Functional Anal. 187, 1 (2001), 200--210. 

\bibitem{ZZZ} J.~Chen, Q.~Leng, X.~Zhao, The Douglas question on the Bergman and Fock spaces, 
arXiv:2405.05412.


\bibitem{CucVas} $\check{\mbox{Z}}$.~$\check{\mbox{C}}$u$\check{\mbox{c}}$kovi\'c, A.~Vasaturo, Carleson measures and Douglas’ question on the Bergman space, Rend. Circ. Mat. Palermo, II. Ser 67 (2018), 323--336. 

\bibitem{D} R.~Douglas, Banach Algebra Techniques in the Theory of Toeplitz Operators, Amer. Math.
Soc., Providence, RI 1980.

\bibitem{GZ} N.~Guan, X.~Zhao, Invertibility of Bergman Toeplitz operators with harmonic polynomial symbols, Sci. China Math. 63 (2020), 965--978.

\bibitem{GZZ} K.~Guo X.~Zhao, D.~Zheng, The spectral picture of Bergman-Toeplitz operators with harmonic polynomial symbols, Ark. Mat 61 (2023), 343--374.

\bibitem{Hal} P.~Halmos, A Hilbert space problem book, 2nd ed. Graduate texts in mathematics
vol. 19, Springer-Verlag, New York-Berlin-Heidelberg, 1982. 

\bibitem{Lue} D.~Luecking, Inequalities on Bergman spaces,  Illinois J. Math. 25, 1 (1981), 1--11.

\bibitem{McS} G.~McDonald, C.~Sundberg, Toeplitz operators on the Disc, Indiana Univ. Math. J. 28, 4 (1979), 
595--611.

\bibitem{Yo} R.~Yoneda, Invertibility of Toeplitz operators on the Bergman spaces with
harmonic symbols, J. Math. Anal. Appl. 516 (2022) 126515. 

\bibitem{ZZ14} X.~Zhao, D.~Zheng, Positivity  of Toeplitz operators via
Berezin transform. J. Math. Anal. Appl. 416 (2014), 881--900. 


\bibitem{ZZ16} X.~Zhao, D.~Zheng, Invertibility of Toeplitz operators via
Berezin transform. J. Operator Th. 75, 2 (2016), 475--495.

\bibitem{ZZ19} X.~Zhao, D.~Zheng, Toeplitz operators and the Berezin transforms, in: K. Zhu (Ed.),
Handbook of Analytic Operator Theory, Chapman and Hall/CRC, New York, 2019.

\bibitem{Zh} K.~Zhu, Operator Theory in Function Spaces, 2nd ed. Mathematical surveys and monographs 
vol. 138, American Mathematical Society, Providence RI, 2007.




\end{thebibliography}
\end{document}